\documentclass{amsart}
\usepackage{latexsym,amsxtra,amscd,ifthen}
\usepackage{amsfonts}
\usepackage{verbatim}
\usepackage{amsmath}
\usepackage{amsthm}
\usepackage{amssymb}

\numberwithin{equation}{section}

\theoremstyle{plain}
\newtheorem{theorem}{Theorem}[section]
\newtheorem{theorema}[theorem]{Theorem}
\newtheorem{theordef}[theorem]{Theorem--Definition}
\newtheorem{prop}[theorem]{Proposition}
\newtheorem{lemma}[theorem]{Lemma}
\newtheorem{cor}[theorem]{Corollary}

\newtheorem{conj}[theorem]{Conjecture}

\theoremstyle{definition}

\newcommand{\N}{\ensuremath{\mathbb N}}

\newcommand{\PP}{\ensuremath{\mathbb P}}
\newcommand{\Z}{\ensuremath{\mathbb Z}}

\def\la{\langle}
\def\ra{\rangle}

\newcommand{\cd}{\mathop{\mathrm{cd}}}

\newcommand{\cohm}{\mathop{\mathrm{cmod}}\nolimits}
\newcommand{\cohp}{\mathop{\mathrm{qgr}}\nolimits}
\def\dim{\mbox{dim\,}}
\def\ext{\mbox{Ext\,}}
\def\gd{\mbox{gl. dim\,}}

\newcommand{\Gr}{\mathop{\mathrm{Gr}}}

\newcommand{\Hom}{\mathop{\mathrm{Hom}}}
\def\H{\mbox{H\,}}
\def\id{\mbox{id\,}}

\newcommand{\pd}{\mathop{\mathrm{pd}}\nolimits}
\newcommand{\proj}{\mathop{\mathrm{proj}}\nolimits}
\newcommand{\Qgr}{\mathop{\mathrm{Qgr}}\nolimits}

\def\rk{\mbox{rank\,}}
\def\sh{\mbox{Sh\,}}

\newcommand{\tor}{\mathop{\mathrm{Tor}}\nolimits}
\newcommand{\Tors}{\mathop{\mathrm{Tors}}\nolimits}
\newcommand{\tors}{\mathop{\mathrm{tors}}\nolimits}

\newcommand{\uext}{\underline{\mathop{\mathrm{Ext}}}}
\newcommand{\uH}{\underline{\mathop{\mathrm{H}}}}
\newcommand{\uhom}{\underline{\mathop{\mathrm{Hom}}}}

\def\tenzA{\mathrel{\mathop{\otimes}\limits_A}}

\begin{document}

\title{Coherent algebras \\
%Two-dimensional Gorenstein algebras  \\
and noncommutative projective lines %${\PP}^1$
%\thanks{Partially
%supported by the grant 05-01-01034 of the Russian Basic Research Foundation.}%
}

\author{Dmitri Piontkovski}

      \address{Department of High Mathematics for Economics,
Myasnitskaya str. 20, State University `Higher School of Economics', Moscow 101990, Russia
% Central Institute of Economics and Mathematics\\
%                      Nakhimovsky prosp. 47, Moscow 117418,  Russia
}
\thanks{Partially
supported by the Dynastia Foundation, by the President of Russian
Federation grant MD-288.2007.1, and by the the Russian Basic
Research Foundation project 05-01-01034}

\email{piont@mccme.ru}



\date{\today}

\begin{abstract}
A well-known conjecture says that every one-relator group is
coherent. We state and partly prove a similar statement for graded
associative algebras. In particular, we show that every Gorenstein
algebra $A$ of global dimension 2 is graded coherent.

  This allows us to define a  noncommutative analogue of
  the projective line
$\PP^1$ as a noncommutative scheme based on the coherent
noncommutative spectrum $\cohp A$ of such an algebra $A$, that is,
the category of coherent $A$-modules modulo the torsion ones. This
category is always abelian $\ext$-finite hereditary with Serre
duality, like the category of coherent sheaves on $\PP^1$. In this
way, we obtain a sequence $\PP^1_n $ ($n\ge 2$) of pairwise
non-isomorphic noncommutative schemes which generalize the scheme
$\PP^1 = \PP^1_2$.
\end{abstract}

\maketitle

%We consider associative algebras over a fixed field $k$.

\section{Introduction}

We will consider \N-graded algebras of the form $A = A_0 \oplus
A_1 \oplus \dots $ over a fixed field $k$. All our algebras will
be assumed to be connected (that is, $A_0 = k$) and finitely
generated.
%So, the word `algebra' below means `connected finitely generated \N-graded algebra'.
All vector spaces and modules are assumed to be
$\Z$-graded, all their elements and maps of them are homogeneous.
%A Hilbert series of such algebra is defined as a formal power series $A(z) = \sum_{i \ge 0} (\dim_k A_i z^i)$.

%We call a quotient algebra $B=A/I$ of an algebra $A$ a
%{\it relative noncommutative complete
%intersection} (RNCI),
%if the ideal $I$ is generated by a strongly free set, see Definition~\ref{def-sf} below.

%\begin{theorem}[Corollary~\ref{cor-rnci-coh}]
%\label{rnci-coh}
%Let $B$ be an RNCI of a graded algebra $A$. If the algebra $B$
%is right Noetherian, then the algebra $A$
%is right graded coherent.
%\end{theorem}

Recall that an algebra (respectively, a group) is called coherent
if every its finitely generated ideal (subgroup) is finitely presented, see
Definition~\ref{def-coher} below.  A well-known conjecture says that every one-relator %discrete
group is coherent~\cite{B}. An analogous statement for graded
algebras also seems to be true.

\begin{conj}
Every graded algebra with a single defining relation is graded coherent.
\end{conj}

(For the definition of coherence,
see~subsection~\ref{ss-def-coher} below.) We will prove this
conjecture provided that the relation is quadratic.

\begin{theorem}[Theorem~\ref{one-quad-rel}]
\label{intro-one-quad-rel}
Every graded algebra defined by a single homogeneous quadratic relation is %right and left
graded coherent.
\end{theorem}

Note that there are non-coherent quadratic algebras with two
relations, for example, the algebras $k\langle x,y,z,t | tz-zy,zx
\rangle$~\cite[Prop.~10]{pi3} or even  $k\langle x,y,z | yz-zy,zx
\rangle$~\cite[Example~2]{pol}.

Recall~\cite{zhang} that a graded algebra $A$ is called regular if
it has finite global dimension (say, $d$) and satisfies the
following Gorenstein property:
$$
\ext^i_A(k_A,k_A) \cong\left\{ \begin{array}{ll} 0, & i\ne d \\
k [l] \mbox{ for some } l\in \Z ,& i =d .
\end{array}  \right.
$$
The most important class of regular algebras is the class of
Artin-Shelter (AS) regular algebras, that is, the ones of polynomial
growth. A well-known conjecture~\cite{AS} claims that all these
algebras are Noetherian.

The following conjecture is due to A. Bondal (unpublished).

\begin{conj}
Every regular
%Gorenstein
algebra
%of finite global dimension
is graded coherent.
\end{conj}

Regular
%Gorenstein
algebras of global dimension 2 have been described in~\cite{zhang}.
All these algebras are one-relator.
If such an algebra is generated in degree one, then it is quadratic, but in general
such algebra is only `generalized quadratic' --- like, for example,
the algebra $k\langle x,y | xy-yx = x^3 \rangle$.

\begin{theorema}[Theorem~\ref{gd2-coh-gen}]
\label{intro-gd2-coh-gen}
Every regular
%Gorenstein
algebra of global dimension two
%of finite global dimension
is graded coherent.
%Let $A$ be a Gorenstein algebra of global dimension 2.
%Then $A$ is right and left coherent.
\end{theorema}

Two abelian categories may naturally be associated to any graded
coherent algebra $A$, that is, the category $\cohm A$ of finitely
presented (=graded coherent)  right graded $A$-modules and its
quotient category $\cohp A = \cohm A / \tors A$ by the category
$\tors A$ of finite-dimensional  modules. This category $\cohp A $
plays a role of
%noncommutative analogue of
projective spectrum for noncommutative coherent
algebras~\cite{pol,BvdB}, in generalization of the well-known
construction (due to Artin and Zhang) of noncommutative schemes in
the Noetherian case~\cite{AZ}. In this approach, a noncommutative
projective scheme is a triple
$$
     ( \cohp A , {\mathcal A}, s ),
$$
where  $A$ is a coherent algebra, noncommutative structural sheaf
${\mathcal A}$ is the the image of $ A \mbox{ in } \cohp A $, and
$s$ is the autoequivalence of $\cohp A$ induced by the  shift of
grading. Some details will be given in the
subsection~\ref{b_schemes}.

 The noncommutative
schemes of (Koszul) Noetherian (AS-)regular algebras of global
dimension $n+1$ are usually considered as noncommutative
generalizations of ${\PP}^{n}$. However, in the case of the
projective line ${\PP}^{1}$,  this Noetherian construction does not
give any more than the standard commutative ${\PP}^{1}$  again. On
the other hand,  there are other Noetherian abelian categories whose
properties are close to the ones of the category of coherent sheaves
on
 ${\PP}^{1}$ (that is, they are hereditary $\ext$-finite with Serre
duality)~\cite{rvdb}, but the  ``coordinate rings'' of the
corresponding noncommutative schemes %to these categories
are far from being
connected graded, in contrast to the coordinate ring $k[x_1, x_2]$
of ${\PP}^{1}$.

Here we will introduce another noncommutative generalization of
${\PP}^{n}$, that is, the noncommutative projective schemes
corresponding to (degree-one generated) coherent regular algebras
of dimension $n+1$. We will show in Proposition~\ref{P^n} that the
corresponding $\cohp A$ is an ext-finite category of cohomological
dimension~$n$, and the algebra $A$ (its coordinate ring) may be
recovered by this category via a suitable "representing functor".
In the case of the projective line ${\PP}^{1}$, we can obtain an
infinite sequence $\{ \PP^1_n \}_{n\ge 2}$ of pairwise
non-isomorphic noncommutative schemes analogous to $ \PP^1 =
\PP^1_2$, where the coordinate ring of each $\PP^1_n$ is a
connected graded 2-dimensional algebra with $n$ generators. The
corresponding categories of coherent sheaves are $\ext$-finite
hereditary and satisfy Serre duality and BGG-correspondence. This
is shown in the following

% In the next proposition, we  consider the schemes associated to Koszul
%(or even degree-one generated) $(n+1)$-dimensional regular algebras
%as {\it non-noetherian} analogues of ${\PP}^{n}$. In our coherent
%case, we obtain

% The first properties of these non-noetherian projective
%lines are the following.

%these categories are standard in Noetherian noncommutative
%projective geometry.
%Since we imagine $\cohp A$ as
%the category of coherent sheaves on a noncommutative projective
%variety, we consider here only degree-one generated algebras,
%because of the analogy with the classical theorem of Serre.

\begin{prop}
\label{intro-prop-nc-proj} Let $A$ be a degree-one generated regular
algebra of global dimension 2 with $n \ge 2$ generators.

(a) The categories $\cohm A$ and $\cohp A$ and  the noncommutative
scheme $\PP^1_n = \PP^1_n (k)$ constructed by $A$ are defined (up to
isomorphisms) by the ground field $k$ and the number $n$ only, and
do not depend on the algebra $A$ itself. All these noncommutative
schemes $\PP^1_n$ are pairwise non-isomorphic, with $\PP^1_2 \cong
\PP^1$.

%(b) The algebra satisfy an analogue of the condition
%$\chi$~\cite{AZ}  for coherent modules, in particular, the shifts of
%$A$ form an ample sequence~\cite{pol} in $\cohp A$.

(b) The category $\cohp A$ is $\ext$-finite hereditary with Serre
duality. If $n \ge 3$, then it is not Noetherian, hence it does not
belong to the classification in~\cite{rvdb}.

(c) The bounded derived category ${\mathcal D}^b(\cohp A)$ is
equivalent to the
%The category $\cohp A$ is bounded derived equivalent to the
category of finite $B$-modules modulo
projectives, where $B$ is a commutative Artinian algebra \\
$k[x_1, \dots , x_n]/(x_ix_j, x_i^2 - x_j^2| i \ne j)$.

(d) The bounded derived category ${\mathcal D}^b(\cohp A)$ is
equivalent to the bounded derived category of the
finite-dimensional representations of the generalized Kronecker
quiver $Q_n$ (that is, quiver with two vertices $v_0, v_1$ and n
arrows $v_0 \to v_1$).

%(d) The category $\cohp A$ is also derived equivalent to the
%category presentations of the   of finite $B$-modules modulo
%projectives, where $B$ is a commutative Artinian algebra \\
%$k[x_1, \dots , x_n]/(x_ix_j, x_i^2 - x_j^2| i \ne j)$.

%
%There are  equivalences of bounded derived categories $D^b (\cohm A)
%\simeq  D^b (\cohm B)$ and $D^b (\cohp A) \simeq  D^b (\cohp B)$,
%where $B$ is a commutative Artinian algebra $k[x_1, \dots ,
%x_n]/(x_ix_j, x_i^2 - x_j^2| i \ne j)$. Here the category $D^b
%(\cohp B)$ is equivalent to the
%%quotient category of the
%category $\underline {\grm } B$
%%$\grm B$
%of finite $B$-modules modulo projectives. %free modules.
\end{prop}

%\begin{probl}
%Describe $\cohp A$.
%\end{probl}

%\begin{rema}
%The coherent regular rings inherit many properties of  Noetherian
%ones. For example, these rings satisfy an analogue of the condition
%$\chi$~\cite{AZ}  for coherent modules (that is, $\dim_k
%\ext^i_k(k,M) < \infty$ for all $i \in \N , M\in \cohm A$); the
%proof is almost verbosely the same as in~\cite[Theorem~8.1]{AZ}. It
%follows, in particular, that the degree shift is  ample in $\cohp
%A$, due to a coherent version of the virtue of the Artin--Zhang
%theorem announced in~\cite[remark after Th.~2.4]{pol}.
%
%However, it is not clear that every object in $\cohp A$ is
%represented by a torsion-free module. Therefore, the following
%question arises: what are the points in noncommutative
%non-Noetherian projective lines $\PP^1_n$?
%\end{rema}

The last statement $(d)$ of this proposition shows that the
category of coherents sheaves on our non-commutative projective
line $\PP^1_n$   is derived equivalent to the category of
coherents sheaves on the noncommutative projective space
$\PP^{n-1}$ in the sense of Kontsevich and Rosenberg,
see~\cite[3.3]{kr}. This statement $(d)$ has been originally
commuticated me by Michel Van den Bergh~\cite{vdb}; another proof
is given recently by Minamoto~\cite{m}.

This paper is organized as follows. In section~\ref{s-back}, we
will give a background on coherent algebras, regular algebras of
global dimension 2, and (relative) noncommutative complete
intersections. In section~\ref{s-crit-coher}, we will give the
following criterion for coherence:  if an algebra $B = A/I$ is a
{\it relative noncommutative complete intersection}  of $A$ (that
is, the ideal $I$ is generated by a strongly free set), and $B$ is
right Noetherian, then $A$ is graded coherent. In the next
section~\ref{s-one-ref-coh}, we will apply the above criterion in
order to prove Theorems~\ref{intro-one-quad-rel}
and~\ref{intro-gd2-coh-gen}. Finally, in section~\ref{s-p1} we
will consider noncommutative schemes associated to coherent
regular algebras. In particular, we will prove
Proposition~\ref{intro-prop-nc-proj}.

\subsection*{Acknowledgement}

I am grateful to MPIM Bonn and IHES for their hospitality during
preparation of this paper. Also, I am grateful to Alexey Bondal,
Alexander Polishchuk, and Michel Van den Bergh for helpful
%conversations.
communications.

\section{Backgroung}

\label{s-back}

\subsection{Coherence}

\label{ss-def-coher}

A finitely generated (f.~g.) right module $M$ is called {\it
coherent} if every its  finitely generated submodule is finitely
presented (that is, presented by a  finite number of generators and
relations). Analogously, a graded f.~g. module is called graded
coherent if every its graded f.~g. submodule is finitely presented.
In fact, this notion had been introduced by Serre~\cite{serre} in a
more general case of coherent sheaves.

\begin{theordef}
\label{def-coher} A (graded) algebra $A$ is called (graded) right
coherent, if the following equivalent conditions hold:

(i) every (homogeneous) finitely generated right-sided ideal in $A$
is finitely presented, that is, $A$ is (graded) coherent as a right
module over itself;

(ii) every finitely presented (graded) right $A$-module is (graded)
coherent;

(iii) all finitely presented (graded) right $A$-modules form an
abelian category.
\end{theordef}

The proof of equivalence may be found in~\cite{chase} (see
also~\cite{faith}). For example, every right Noetherian algebra is
right coherent, as well as every free associative algebra and
every path algebra.
%, and finitely presented monomial algebras are coherent as
% well~\cite{pi1}.

Since all our algebras and modules  are  graded, by the word {\it
coherent} we will mean {\it graded  right coherent} algebras  and
modules. The idea of noncommutative geometry based on such
algebras will be explained in the next subsection.

\subsection{Noncommutative schemes}

\label{b_schemes}

 Let $A$ be a graded algebra. By $\Gr A$ (respectively, $\cohm
A $) we denote the abelian category of graded (resp., coherent)
$A$-modules. Let $\Tors A$ (resp., $\tors A$) be the category of
torsion $A$-modules (resp., finite dimensional modules), where a
module $M$ is called torsion if for every $x \in M$ there is $n>0$
such that $x A_{\ge n} =0$. Note that $\Tors A$  is a Serre
subcategory of $\Gr A$; moreover, if $A$ is coherent, then $\tors A$
is also a Serre subcategory of $\cohm A$. The quotient abelian
categories $\Qgr A = \Gr A / \Tors A$ and $\cohp A = \cohm A / \tors
A$ (for
 coherent $A$) play the roles  of
 the categories %$\Proj A$ and $\proj A$
 of
 (quasi)coherent sheaves on the projective scheme  associated to $A$.
Due to classical Serre theorem~\cite{serre}, these categories of
modules are indeed equivalent to the respective categories of
sheaves provided that $A$ is commutative.

The image ${\mathcal A}$ of $A_A$ in  $\Qgr A$ (or in $\cohp A$)
plays the role of the structure sheaf, and the the degree shift $s:
M \mapsto M[1]$ plays the role of the polarization. Thus, a %general
noncommutative scheme is a triple
%Let $X$ be a  triple
 $$
X = ({\mathcal C}, {\mathcal A}, s ),
 $$
 where ${\mathcal C}$ is a suitable
 $k$-linear abelian category, ${\mathcal A}$ is an object, and $s$ is an
autoequivalence of ${\mathcal C}$. For ${\mathcal C} = \Qgr A$ with
an arbitrary connected graded algebra $A$, this definition is due to
Verevkin~\cite{ver} (a general scheme). For ${\mathcal C} = \cohp A$
(coherent scheme), this definition is due to Artin and
Zhang~\cite{AZ} in the case of noetherian $A$ (noetherian scheme)
and to Bondal and Van den Bergh~\cite{BvdB} and
Polishchuk~\cite{pol} in a more general setting of coherent algebra
$A$.

According to~\cite{AZ}, a morphism $f: X \to X'$ of two schemes $X =
({\mathcal C}, {\mathcal A}, s )$ and $X' = ({\mathcal C'},
{\mathcal A'}, s' )$ is a $k$-linear functor $f: {\mathcal C} \to
{\mathcal C'}$ such that $f({\mathcal A})$ is isomorphic to
${\mathcal A'}$ and there is an isomorphism of functors $fs \cong s'
f$. A map of schemes is defined as an isomorphism class of
morphisms. Such a morphism $f$ (or a respective map) is called an
isomorphism if it is an equivalence of categories $f: {\mathcal C}
\cong {\mathcal C'}$.

Given such a triple $X = ({\mathcal C}, {\mathcal A}, s )$, we can
apply an analogue of the Serre functor to define a connected graded
algebra $A:= \Gamma_{\ge 0} (X) = \bigoplus_{i\ge 0} \Hom ({\mathcal
A}, s^i ({\mathcal A}))$ with the multiplication $ a \cdot b :=
s^j(a) \circ b$ for $a : {\mathcal A}\to s^i ({\mathcal A}), b:
{\mathcal A}\to s^j ({\mathcal A})$. In some cases, this  algebra
$A$ is coherent and the scheme $X$ itself is isomorphic to the
scheme $(\cohp A, {\mathcal A}, s )$. This happens if the
autoequivalence $s$ is ample~\cite{AZ}, that is, the shifts of
${\mathcal A}$ form an ample sequence in ${\mathcal C}$~\cite{pol}.

If two general schemes $X$ and $Y$ are isomorphic, then   the
algebra $\Gamma_{\ge 0} (Y)$ is isomorphic to a Zhang twist of
$\Gamma_{\ge 0} (X)$; on the other hand, if a coherent algebra $B$
is a Zhang twist of an algebra $A$, then the coherent (and general)
schemes of these algebras are isomorphic~\cite[Th.~1.4]{zh-tw}. Here
an algebra $B$ is called a Zhang twist of $A$ if there are
$k$-linear bijections $\tau_i: A_i \to B_i , i \ge 0$ such that
$\tau_{m+n} (y z)= \tau_m(y) \tau_{m+n}(z)$ for homogeneous $y\in
y\in A_n, z\in A$~\cite[Prop.~2.8]{zh-tw}. For example, the
projective scheme of the quantum polynomial algebra $k \la x,y | xy
= q yx \ra $ is isomorphic to $\PP^1$ for every $q \ne 0$.

Let $A$ be a graded algebra, let $M,N \in \Gr A$ be two modules, and
let ${\mathcal M}$ and ${\mathcal N}$ be their images in $\Qgr A$.
Let $\uhom ({\mathcal M},{\mathcal N}) := \bigoplus_{i\in \Z} \Hom
 ({\mathcal M}, s^i{\mathcal N})$,
and let $\ext$ and $\underline \ext$ be the derived functors of
$\Hom$ and $\uhom$. Since the obvious functor $\cohp A \to \Qgr A$
is fully faithful for a coherent algebra $A$, the functors $\ext$
and $\uext$ on the category $\cohp A$ are restrictions of the
respective functors on $\Qgr A$.

Following~\cite{ver, AZ}, we define the cohomologies of objects
of $\Qgr A $ as $\H^i ({\mathcal M})= \ext^i ({\mathcal A},
{\mathcal M}) $ and $\uH^i ({\mathcal M})= \underline \ext^i
({\mathcal A}, {\mathcal M}) = \lim\limits_{n \to \infty} \ext^i_A
(A_{\ge n}, M) $.  According to~\cite[Lemma~4.1.6]{BvdB}, we have
$\uH^i ({\mathcal M}) \cong \lim\limits_{n \to \infty} \ext^i_A
(A_{\ge n}, M) $ and $\H^i ({\mathcal M}) \cong \lim\limits_{n \to
\infty} \ext^i_A (A_{\ge n}, M)_0$.

\subsection{Regular algebras of global dimension 2}

Let us recall some results of~\cite{zhang}. Let $V $ be a vector
space. A {\it rank} of an element $b \in T(V)$  is defined as the
minimal number $r$ of elements $l_1, \dots, l_r \in V$ such that $b
= l_1 a_1 +\dots + l_r a_r$ for some $a_1 , \dots , a_r \in T(V)$.

\begin{theorema}[\cite{zhang}]
 \label{zh-gen}
A graded algebra $A$ is regular of global dimension 2 if and only if it is isomorphic to the algebra
$k\langle x_1 , \dots , x_n \rangle / (b)$, where $\rk b =n >1$, or, equivalently, the following
conditions hold:

1.  $n\ge 2$;

2. $1 \le \deg x_1 \le \dots \le \deg x_n$ with $\deg b = \deg x_i + \deg x_{n+1-i}$
for all $i$;

3.  for some graded automorphism $\sigma$ of the free algebra
$k\langle x_1 , \dots , x_n \rangle $ we have $b = \sum_{i= 1}^n x_i
\sigma(x_{n+1-i}) $.

In this case, the algebra $A$ is Noetherian if and only if $n=2$.
\end{theorema}

%See the proof in \cite[Theorems~0.1,~0.2 and Prop.~1.1]{zhang}.

In particular, ir follows that a regular two-dimensional algebra is
Koszul if and only if it is degree-one generated.

\subsection{(Relative) noncommutative complete intersections}

In the next definition, we will unite several statements
from~\cite{an1}. For discussions on strongly free sets as a
noncommutative analogue of regular sequences and related topics,
see also~\cite{pi2, ufn}. Recall that a relative complete
intersection, from an algebraic point of view, is a quotient of
some graded or local commutative ring by an ideal generated by a
regular sequence. Here we
 introduce a
 % the notion of
 relative noncommutative complete
intersection (RNCI) as a quotient of a graded algebra by an ideal
generated by a strongly free set. It is analogous to the term
`noncommutative complete intersection', that is, RNCI of a free
algebra~\cite{an1, golod, EG}.

%For the discussion on Shafarevich complex, see~\cite{pi2}.

\begin{theordef}[\cite{an1}]
\label{def-sf} Suppose that a set $X$ of homogeneous elements in a
graded algebra $A$ minimally generates a two-sided ideal $I$. Let
$B=A/I$ be a quotient algebra. The set $X$ is called strongly free,
if
%either of
the following equivalent conditions hold:
%are equthere are isomorphisms of vector spaces
% natural morphism $A\to A/I$
%induces isomorphisms of graded vector spaces

(i) there are isomorphisms of graded vector spaces
$$
\tor_i^A(k,k) \simeq \tor_i^{B}(k,k) \mbox{ for all $i \ge 3$ and }
$$
$$
    \tor_1^A(k,k) \oplus \tor_2^{B}(k,k) \simeq  \tor_1^B(k,k) \oplus \tor_2^{A}(k,k) \oplus kX;
$$

(ii) there is an isomorphism of graded vector spaces
$$
      B\langle X \rangle \simeq A,
$$
where $ B\langle X \rangle  = B * k\langle X \rangle$ is a free product of $B$ and  a free algebra on $X$;

(iii)
the Shafarevich complex $\sh (X, A)$ is acyclic in positive degrees;

%(iv) the following equality of formal power series holds:
%$$
%     A(z)^{-1} = B(z)^{-1} + X(z),
%$$
%where $X(z) = \sum\limits_{x \in X} z^{\deg x}$.

In this situation, we refer to the algebra $B=A/I$  as a relative
noncommutative complete intersection (RNCI) of the algebra $A$.
\end{theordef}

In particular, it follows that if $B=A/I$ and there are isomorphisms
$\tor_i^A(k,k) \simeq \tor_i^{B}(k,k)$ for all $i \ge 2$, then $B$
is an RNCI of $A$.

%Note that the property~$(ii)$ here may be re-written as follows. For
%every there is an isomorphism of graded vector spaces $\alpha :
%B\langle X \rangle \to A$ such that $\alpha (B) = B'$, $\alpha ( B X
%B\langle X \rangle ) = I$.

\section{A criterion for coherence}

\label{s-crit-coher}

\begin{lemma}
Let $X$ be a strongly free set in a graded algebra $A$ and  let $I$
be an ideal generated by $X$. Then $I$ is  free as right (and left)
$A$-module.

More precisely, let  $B = A/I$, and
let $B'$ be any pre-image of $B$ in $A$ with the natural isomorphism of vector spaces $B' \cong B$.
Then $I$ as a free right $A$-module is minimally generated by the vector space $B'X$.
\end{lemma}

\begin{proof}
Obviously, $I = B'X A$. We have to show that the natural
epimorphism $\gamma : B' X \otimes_k A \to B'X A$ is an
isomorphism. Following~\cite{an1}, there is an isomorphism of
graded vector spaces $\alpha : B\langle X \rangle \to A$ such that
$\alpha (B) = B'$,  $\alpha ( B X B\langle X \rangle  ) = I$ (this
follows also from the property~$(ii)$ of Theorem~\ref{def-sf}).
The right $B\langle X \rangle $-module in the last equality is
free, so, we get the desired isomorphisms of graded vector spaces
$ B'X A = I \cong B X B\langle X \rangle = B' X \otimes B\langle X
\rangle \cong B' X \otimes_k A$. Therefore, the map $\gamma $ is
an isomorphism of graded vector spaces. So, it is an isomorphism
of modules.
\end{proof}

The next statement is similar to \cite[Prop.~3.3]{pol}.

\begin{prop}
\label{prop-coh-crit}
 Let $B = A/I$, where the algebra $B$ is right
Noetherian and the ideal $I$ is  free as a left $A$-module. Then the
algebra $A$ is right graded coherent.
\end{prop}

\begin{proof}
Let $J$ be a proper finitely generated homogeneous right-sided ideal in $A$.
% generated by a finite homogeneous set $Y$.
We have to show that $J$ is finitely presented, that is, $\dim_k \tor^A_2(A/J, k) < \infty$.

Consider a standard spectral sequence $E^2_{p,q}
=\tor_p^B(\tor_q^A(A/J,B),k) \Longrightarrow \tor_* (A/J,k)$. Let
$N_q = \tor_q^A(A/J,B)$. Since the projective dimension of the
left module ${}_A\!B$ is at most one (because it admits a free
resolution $0\to I \to A$), we have $N_q= 0 $ for $q>1$, hence
$E^2_{p,q} =0$ for $q>1$. The right $B$-module $N_0 =  A/J
\otimes_A B$ is obviously finitely generated. Moreover, the short
exact sequence
$$
   0 \to J \to A \to A/J \to 0
$$
gives, after tensoring by B,  an exact sequence
$$
0 \to N_1 \to J \tenzA B \to A \tenzA B \to N_0 \to 0.
$$
Since $N_1$ is a submodule of a finitely generated $B$-module $J \tenzA B$,
it is finitely generated as well. Therefore, we have
$\dim_k E^2_{p,q} = \dim_k \tor_p^B(N_q,k) < \infty$ for all $p, q$.
Thus, $\dim_k \tor^A_2(A/J, k) \le \dim_k E^2_{2,0} +\dim_k E^2_{1,1} <  \infty$.
\end{proof}

\begin{cor}
\label{cor-sf_coher}
\label{cor-rnci-coh} Let $B$ be an RNCI of a
graded algebra $A$. If the algebra $B$ is right Noetherian, then the
algebra $A$ is right graded coherent.
%Let $X$ be a strongly free set in an algebra $A$ and  let $I$ be an ideal generated by $X$.
%If the quotient algebra $B = A/I$ is right Noetherian, then $A$ is right graded coherent.
\end{cor}

%\begin{rema}
%In fact, it may be shown more:  in the both
%Proposition~\ref{prop-coh-crit} and Corollary~{cor-sf_coher}, the
%algebra $A$  is also coherent in the non-graded sense.
%\begin{rema}

\section{One-relator quadratic algebras}

\label{s-one-ref-coh}

\begin{theorem}
\label{one-quad-rel}
Every algebra defined by a single homogeneous quadratic relation is %right and left
graded coherent.
\end{theorem}

\begin{proof}
Let $b \in V \otimes V$ be a quadratic element in the free algebra
$T(V)$, where $\dim V = n$, and let $A$ be a quotient algebra of
$T(V)$ by an ideal ${\id (b)}$ generated by $b$. If $n=1$, then
the algebra $A$ is finite-dimensional, hence Artinian, hence
Noetherian, and coherent. If $n=2$, then either $A$ is Noetherian
or $b$ has the form $b=xy$, where $x,y \in V$~\cite[p.~172]{AS}.
In the last case, $A$ is coherent by~\cite[Th.~2]{pi1}.

Consider the case $n\ge 3$.
%Suppose that  $b = \sum_{i=1}^r x_i y_i$, where $x_i,y_i \in V$.
%The minimal such $r$ is called the {\it rank} of $b$~\cite{zhang}.
If $\rk b =1$, then $b = xy$ is a monomial on generators, and the
algebra $A$ is coherent, again by~\cite[Th.~2]{pi1}. So, we can
assume that $\rk b \ge 2$.

\begin{lemma}
\label{quadr-decomp}
Let $V$ be an $n$-dimensional vector space with $n \ge 2$, and let $b \in V \otimes V$.
Given an $(n-2)$-dimensional subspace $W\subset V$, let $b'$ be the image of $b$ in
$(V/W) \otimes (V/W)$. Then either $b = xy$ for some $x,y \in V$ or there exists $W$
such that $\rk b' = 2$.
%In the last case, there is a basis $\{ a,b \}$ of the vector space  $V/W$
%such that $b' = ab + \alpha ba + \beta b^2$ for some $\alpha , \beta \in k, \alpha \ne 0$.
\end{lemma}

\begin{proof}[Proof of Lemma~\ref{quadr-decomp}]
%Let $x_1 , \dots , x_n$ be a basis of $V$.
%, and let $b = \sum_{i,j} \alpha_{ij} x_i x_j$.
By the induction on $n$, we can assume that for every $x \in V$,
the image $b''$ of $b$ in $(V/kx) \otimes (V/kx)$ has rank $\le
1$. Let $\{ x_1=x , \dots , x_n \}$ be a basis of $V$.

If, for some $x$, we have $b'' =0$, then $b = \alpha x^2 + x l_1 + l_2 x$ with
$l_i \in k\{x_2 , \dots , x_n \}, \alpha \in k$. If $\rk b \ge 2$, then
 $l_1 \ne 0 $ and $l_2 \ne 0$. Now, if $l_2 \ne \beta l_1 $ for $\beta \in k$,
 then the image $b = \alpha x^2 + x l_1 + l_1 x$ has rank 2 --- a contradiction.
 Hence, $l_2 = \beta l_1 $ for some  $0 \ne \beta \in k$,
 and $b = \alpha x^2 + x l_1 + \beta l_1 x$ has rank two.

  So, we can assume that $\rk b'' =1$ for every $x$.
  Then $b'' = uv$ for some nonzero $u,v \in k\{x_2 , \dots , x_n
  \}$. Hence $b = uv + \alpha x^2 + x l_1 + l_2 x$ with $l_i \in k\{x_2 , \dots , x_n \}, \alpha \in k$.
  We can assume that either (1) $u=x_2, v=x_3$ or (2) $u=v=x_2$.

  Suppose that $l_1 \ne 0$ and $l_2 \ne 0$.
   The image of $b$ under the factorization by $l_1$ has unit rank,
   hence $l_2 = \beta u$ for some $\beta \in k$.
  Similarly, $l_1 = \gamma v$ with $\gamma \in k$.
   In the case (1), let $W = k\{ x_2-x_3, x_4, \dots , x_n\}$;
   in the case (2), let us put $W = k\{ x_3, \dots , x_n\}$.
   In both cases, the image $b'$ of $b$ in $(V/W) \otimes (V/W)$ has the same rank as $b$.

   Now, it remains to consider the case $l_2 = 0$ (the case $l_1 = 0$ is analogous).
   Then $b = uv + x( \alpha x + l_1)$. If $\alpha = 0$ and $l_1=0$, then $\rk b=1$, and there is nothing to prove.
   In the case (1), the image of $b$  under the factorization by $(x_2-x_3)$
   must have rank one, hence $\alpha =0$, $l_1 = \lambda v$ for some $\lambda \in k$, and $\rk b =1$.
   In the case (2), because either $l_1=0$ or the the image of $b$
    under the factorization by $l_1$ has unit rank, we have  $l_1 = \lambda x_2$
     for some $\lambda \in k$. Then $b$ depends on the variables $x_1$ and $x_2$ only, hence
   we may put $W = k\{ x_3,  \dots , x_n\}$.
  \end{proof} %\begin{proof}[Proof of Lemma~\ref{quadr-decomp}]

Recall that $\rk b \ge 2$.
Let $x_1 , \dots , x_n$ be a basis of $V$ such that $W = k\{ x_i | i =3 \dots n \}$
be as in this Lemma. Then the image $b'$ of $b$ in $(V/W) \otimes (V/W) = k\{ x_1, x_2\}^{\otimes 2}$
has rank 2.
By Theorem~\ref{zh-gen}, the algebra $B = A/ \id (x_3, \dots, x_n) =
k\langle x_1, x_2 | b' =0 \rangle$ is  Noetherian. % domain.
Now, the set $X = \{ x_3, \dots, x_n \}$ is strongly free in the algebra $A $,
because $A/\id (X) = B$, while $\tor_1^B(k,k) \oplus kX  \cong \tor_1^A(k,k)$,
 $\tor_2^B(k,k) \cong \tor_2^A(k,k) \cong kb$, and $\tor_i^B(k,k) = \tor_i^A(k,k) =0$ for all $i\ge 3$.
Thus, it follows from Corollary~\ref{cor-sf_coher} that  the algebra
$A$ is coherent.
\end{proof}

\begin{theorema}
\label{gd2-coh-gen} Let $A$ be a  regular algebra of global
dimension 2.
%Gorenstein algebra of global dimension 2.
Then $A$ is graded coherent.
\end{theorema}

\begin{proof}
According to Zhang's Theorem~\ref{zh-gen}, the algebra $A$ has the
form $A = k\langle x_1 , \dots , x_n \rangle / (b)$, where $n\ge
2$, $1 \le \deg x_1 \le \dots \le \deg x_n$ with $\deg b = \deg
x_i + \deg x_{n+1-i}$ for all $i$, and for some graded
authomorphism $\sigma$ of the free algebra $k\langle x_1 , \dots ,
x_n \rangle $ we have $b = \sum_{i= 1}^n x_i \sigma(x_{n-i}) $. If
$\deg  x_1 = \dots = \deg x_n$, then $b \in k\{ x_1 , \dots , x_n
\}^{\otimes 2}$. Hence $A$ is coherent by
Theorem~\ref{one-quad-rel}.

So, we can assume that $\deg x_1 = \dots = \deg x_p < \dots < \deg
x_{n-p+1} =\dots =\deg x_n$. Since the definition of regular rings
is left-right symmetric, it follows that there is another graded
automorphism $\tau$ of the free algebra $k\langle x_1 , \dots ,
x_n \rangle $ such that $b = \sum_{i= 1}^n \tau (x_{n-i}) x_i $.
Let $\tilde b = \sum_{i=n-p+1}^n (x_i \sigma(x_{n-i}) + \tau
(x_{n-i}) x_i)$. Obviously, the element $b-\tilde b$ does not
depend on the variables $x_{n-p+1}, \dots ,x_n$, hence $\rk
(b-\tilde b) \le n-p$ (where rank is defined as the minimal number
$r$ of elements $l_1, \dots, l_r \in k\{ x_1 , \dots , x_n \}$
such that $b-\tilde b = l_1 a_1 +\dots l_r a_r$ for some $a_1 ,
\dots , a_r \in k\langle x_1 , \dots , x_n \rangle$). Now, we are
interested in $\rk \tilde b$.

Consider the case  where $\rk \tilde b \le 1$. Then $\rk b \le \rk
(b-\tilde b) + \rk \tilde b \le n-p+1$. Since $\rk b =n$
by~Theorem~\ref{zh-gen}, we have $p=1$. Then $\tilde b  b = x_n
\sigma(x_1) + \tau (x_1) x_n = \alpha x_n x_1 +\beta x_1 x_n$ for
some nonzero $\alpha, \beta \in k$. Thus, $\rk \tilde b =2$ --- a
contradiction.

So, $\rk \tilde b \ge 2$. Note that $\tilde b \in V \otimes V$, where
$V = k\{ x_1 , \dots , x_p, x_{n-p+1}, \dots , x_n \}$. According to Lemma~\ref{quadr-decomp},
there is a $(2p-2)$-dimensional  subset   $W$ in $V$
(say, $W = k\{ x_2 , \dots , x_p, x_{n-p+1}, \dots , x_{n-1} \}$)
such that the rank of the image $b'$ of $\tilde b$ in  $ (V/W) \otimes (V/W)$
is 2. It follows from~Theorem~\ref{zh-gen}
that the algebra $B = k\langle x_1, x_n | b' \rangle = A/\id (x_2, \dots, x_{n-1})$
is Noetherian and has  global dimension 2.
By the same arguments as in the proof of Theorem~\ref{one-quad-rel}, %we conclude that
the set $X = \{ x_2, \dots, x_{n-1} \}$ is strongly free in $A$. In
the view of Corollary~\ref{cor-sf_coher}, we conclude that the
algebra $A$ is coherent.
\end{proof}

\section{Non-noetherian ${\PP}^1$}

\label{s-p1}

A module $M$ over an algebra $R$ is said to satisfy condition $\chi$
if $\dim_k \ext^i (k,M)<\infty$ for all $i \ge 0$, see~\cite{AZ}. A
coherent algebra $A$ said to satisfy $\chi$ if every finitely
presented $A$-module $M$ satisfy $\chi$.

The following proposition is similar to~\cite[Th.~8.1]{AZ}. The
proof is more or less similar too.

\begin{prop}
\label{P^n} Let $A$ be a graded coherent regular algebra of global
dimension $d \ge 0$. Then

(1) $A$ satisfies the condition $\chi$;

(2) the algebra $A$ may be recovered from its noncommutative scheme
$\proj A := ( \cohp A , {\mathcal A}, s )$ as
$$
A \cong \Gamma_{\ge 0} (\proj A);
%\bigoplus_{i\ge 0} \Hom ({\mathcal A}, s^i ({\mathcal A}));
$$

(3) the category $\cohp A$ is $\ext$-finite and has cohomological
dimension $d-1$.
\end{prop}

Notice that the condition~$(2)$  here means that  $s$ is
ample~\cite{AZ}, that is, that the shifts of ${\mathcal A}$ form an
ample sequence in  $\cohp A$~\cite{pol}.

\begin{proof}
% This proof is similar to the one in the Noetherian case,
% see~\cite[Th.~8.1]{AZ}. It is based also on an analogue
% of~\cite[Prop.~7.2]{AZ} for non-noetherian
% algebras~\cite[Section~4]{BvdB}.
%  the ideas of Bondal and Van den Bergh

Using  the induction on the projective dimension $p$ of a coherent
module $M$, we will show that $\dim_k \ext^i (k,M)<\infty$ for all
$i \ge 0$. If $p=0$, then $M$ is a finitely generated free
$A$-module, so, all $\ext ^i (k,M)$ are bounded by the Gorenstein
condition. If $p>0$, then there is a short exact sequence
(presentation)
%$$
\begin{equation} 0\to N \to P \to M \to 0
\label{simple_short_exact}
\end{equation}
%$$
where $P$ is a projective module and $\gd N <p$. By the induction
assumption, the condition $\chi$ holds for $P$ and $N$; by the
exact triangle of $\ext$s, it holds for $M$ as well.

$(2)$. Let ${\mathcal M}$ be an image of some $M \in \cohm A$ in
$\cohp A$.
%Let us denote $H^i ({\mathcal M})= \ext^i ({\mathcal A}, {\mathcal
%M}) = \lim\limits_{n \to \infty} \ext^i_A (A_{\ge n}, M)_0 $ and
%$\underline H^i ({\mathcal M})=  \underline \ext^i ({\mathcal A},
%{\mathcal M}) = \lim\limits_{n \to \infty} \ext^i_A (A_{\ge n}, M) $
%(the rightmost isomorphisms follow from~\cite[Lemma~4.1.6]{BvdB}).
For every $n>0$, the short exact sequence $0\to A_{\ge n} \to A \to
A/A_{\ge n} \to 0$ gives an exact sequence
$$
   0\to \Hom_A(A/A_{\ge n}, A) \to A\to \Hom_A(A_{\ge n}, A) \to \ext^1_A (A/A_{\ge n},
             A)\to 0.
$$
Since $A$ is regular, the left and right terms are zero. Hence $A
\cong \Hom_A(A_{\ge n},A) \cong \lim\limits_{n \to \infty} \Hom_A
(A_{\ge n}, A) = \uH^0 ({\mathcal A})$.

$(3)$.
 Notice that  $\uH^i({\mathcal A}) = \lim\limits_{n \to \infty}
\ext^{i+1}_A (A_{\ge n}, A) $  for $i\ge 1$, hence $\uH^i({\mathcal
A}) = 0$ for $i \ne 0, d-1$ and $\uH^{d-1}({\mathcal A}) = A^*[l]$
for some $l \in \Z$. It follows that $\cd (\cohp A) \ge d-1$ and
that the cohomologies $\H^i({\mathcal A})$ are finite-dimensional.

Moreover, $\uH^i({\mathcal M}) = \lim\limits_{n \to \infty} \ext^i_A
(A_{\ge n}, M)=0$ for all $i \ge d$. If $\pd M =0 $, then $M =
\bigoplus_i A[l_i] $ is a finitely generated free module, hence
$\ext^i ({\mathcal A}[l], {\mathcal M}) $ is finite-dimensional for
every $i\ge 0, l\in \Z$. By the induction on $\pd M$, it follows
from the $\ext ({\mathcal A}[l], -)$ triangle for the exact
sequence~(\ref{simple_short_exact}) that the vector spaces
 $\ext^i ({\mathcal A}[l], {\mathcal M}) $  are
finite-dimensional for all $i, l, M$.

Let ${\mathcal M'}$ be an image in $\cohp A$ of another coherent
$A$-module $M'$. If $M' = \bigoplus_{i=0}^t A[l_i] $ is a free
module, we can apply the functor $\ext^i(-, {\mathcal M}) $ to the
short exact sequence $0\to \bigoplus_{i=0}^{t-1} A[l_i] \to
{\mathcal M'} \to {\mathcal A}[l] \to 0$.  The derived exact
triangle shows that  the vector space $\ext^i ({\mathcal M'},
{\mathcal M}) $ is finite-dimensional for every $i$ and vanishes for
$i\ge d$. For non-free modules $M'$, we
 proceed by induction on $\pd M'$. Applying the same functor to the
short exact sequence $0\to N \to F \to M' \to 0$ analogous
to~(\ref{simple_short_exact}), we deduce  that the vector spaces
$\ext^i({\mathcal M'}, {\mathcal M})$ are finite-dimensional for all
$i$ and vanish for $i \ge d$ as well. It follows that the
cohomological dimension of $\cohp A$ is $d-1$ and that the category
$\cohp A$ is $\ext$-finite.
\end{proof}

\begin{proof}[Proof of Proposition~\ref{intro-prop-nc-proj}]

 According to~\cite[Th.~1.4]{zh-tw}
 (see also the subsection~\ref{b_schemes} above),
the coherent scheme of $A$ is independent  (up to isomorphism) on
the choice of the automorphism $\sigma$ in Theorem~\ref{zh-gen}. On
the other hand, if the regular algebras   $A$ and $A'$ of global
dimension two have different numbers of generators (say, $m$ and
$n$), then they are not twists of each other because $\tau_1: A_1
\to A_1'$ cannot be an isomorphism of vector space. This
proves~$(a)$.

Let us also give  a direct proof of the last statement. Let be
${\mathcal A}$ and ${\mathcal A'}$ the images of $A$ and $A'$ in
respective $\cohp$, and let $s$ and $s'$ be the shifts if grading
in these $\cohp$. Assume that the schemes $\PP^1_n$ and $\PP^1_m$
with the underlying algebras  $A$ and $A'$ are isomorphic.
% for some $m>n$.
By definition~\cite{AZ}, this means that there is an equivalence of
categories $F : \cohp A \to \cohp A'$ such that $F({\mathcal A})
\cong {\mathcal A}'$ and $s'F \cong F s$. Then $F$ maps the exact
sequence
$$
0 \to s^2 {\mathcal A} \to s {\mathcal A}^n \to {\mathcal A} \to 0
$$
to the exact sequence
$$
0 \to s'^2 {\mathcal A'} \to s' {\mathcal A'}^n \to {\mathcal A'}
\to 0.
$$
Taking the Euler characteristics for the second exact sequence, we
deduce that the following equality of formal power series holds for
some polynomial $p(z) \in \Z [z]$ (because a pre-image of this
sequence in $\cohm A'$ must be exact up to finite-dimensional
modules):
$$
         A'(z) (1-n z + z^2) = p(z),
$$
where $A'(z):= \sum_{i\ge 0} (\dim A'_i)z^i = (1-mz+z^2)^{-1} $. It
follows that $m=n$.

%The proofs of the condition $\chi$ for coherent modules over $A$
%(that is, $\dim_k \ext^i_k(k,M) < \infty$ for all $i \in \N , M\in
%\cohm A$)  and of the hereditary property  and Serre duality for
%$\cohp A$ is verbosely the same as the proofs in the Noetherian
%case~\cite[Theorem~8.1]{AZ}. Then the ampleness in~$(b)$ follows
%from a coherent analogue of the virtue of the Artin--Zhang theorem
%announced in~\cite[remark after Th.~2.4]{pol}.

$(b)$
The Serre duality for   $\cohp A$ follows from~\cite{s-dual}.%[Theorem~4.9]
The hereditarity  (that is, that  $\cohp A$ has cohomological
dimension $\le 1$) and $\ext$-finiteness follows from
Proposition~\ref{P^n}.
%Then the ampleness in~$(b)$ follows

If $n\ge 3$, then $A$ is not Noetherian by Theorem~\ref{zh-gen}. Let
us show that the image ${\mathcal A}$ of $A$  in $\cohp A$  is not
Noetherian  as well. In the view of~$(a)$, we may assume that $b =
x_1 x_2 + x_2 x_3 +\dots + x_n x_1$. Then $b$ forms a Groebner basis
of the ideal $\id (b) \subset k\langle x_1 , \dots , x_n \rangle$
w.~r.~t. an arbitrary deg-lex order, therefore, there is a linear
basis of $A$ consisting of the monomials on the variables $x_1,
\dots , x_n$ which do not contain a subword $x_1 x_2$. Now, it is
easy to see that the monomials $x_1 x_3, x_1^2 x_3 , \dots, x_1^t
x_3$ form a right Groebner basis of the right-sided ideal $I_t
\subset A$ generated by them. It follows that every quotient module
$I_{t}/I_{t-1}$ is infinite-dimensional (because it contains a
sequence of linearly independent monomials $x_1^{t}x_3^s, s \ge 1$),
hence the image in ${\mathcal A}$ of the chain $ I_1 \subset I_2
\subset \dots $  is strictly ascending. This proves~$(b)$.

The statement~$(a)$ allows us to choose any particular $b$ of rank
$n$; let us choose $b =x_1^2 + \dots + x^2_n$. Then we have $B =A^!$
--- a Koszul dual algebra. Now, the claim~$(c)$ follows from
the Koszul duality and a noncommutative analogue of the
Bernstein-Gelfand-Gelfand correspondence,
see~\cite[Prop.~4.1 and Cor.~4.5]{mvs}.% and Serre duality.

The claim~$(d)$ has been communicated me by Michel Van den
Bergh~\cite{vdb}. In the view of~$(a)$, it is sufficient to show
the derived  equivalence for the same value $x_1^2 + \dots +
x^2_n$ of $b$. In this case, this equivalence has been also shown
in~\cite[Theorem~0.1]{m}.
\end{proof}

\end{document}